\documentclass[12pt]{amsart}

\usepackage{cite}

\usepackage{amsfonts}
\usepackage{amssymb, latexsym, amsthm, amsmath, verbatim}
\usepackage{mathrsfs}
\usepackage{multirow,makecell,rotating,diagbox}
\usepackage{indentfirst}
\usepackage{float}

\numberwithin{equation}{section}
\theoremstyle{definition}
\newtheorem{Thm}{Theorem}
\newtheorem{Prop}[equation]{Proposition}
\newtheorem{Cor}[equation]{Corollary}
\newtheorem{Lem}[equation]{Lemma}

\newtheorem{Exa}[equation]{Example}

\makeatletter
\def\imod#1{\allowbreak\mkern5mu{\operator@font mod}\,\,#1}
\makeatother

\usepackage{color}

\definecolor{blue}{rgb}{0,0,1}
\definecolor{red}{rgb}{1,0,0}
\definecolor{green}{rgb}{0,.6,.2}
\definecolor{purple}{rgb}{1,0,1}

\long\def\red#1\endred{{\color{red}#1}}
\long\def\blue#1\endblue{{\color{blue}#1}}
\long\def\purple#1\endpurple{{\color{purple}#1}}
\long\def\green#1\endgreen{{\color{green}#1}}

\begin{document}

\title[]{Dimension formulas of modular form spaces with character for Fricke groups}
\author{Yichao Zhang$^\star$ and Yang Zhou}
\address{School of Mathematics, Harbin Institute of Technology, Harbin 150001, P. R. China}
\email{yichao.zhang@hit.edu.cn}
\address{School of Mathematics, Harbin Institute of Technology, Harbin 150001, P. R. China}
\email{18B912038@stu.hit.edu.cn}
\date{}


\begin{abstract} We first obtain the dimension formulas for the spaces of holomorphic modular forms  with character for the Fricke group $\Gamma_0^+(N)$, then that for $\Gamma_0^*(N)$ with all Atkin-Lehner involutions added in a particular case.
\end{abstract}
\maketitle

%\tableofcontents
\newcommand{\Z}{{\mathbb Z}} % for integers
\newcommand{\Q}{{\mathbb Q}} % for rational
\newcommand{\R}{{\mathbb R}}
\newcommand{\m}{{\text{ mod }}}

\section{Introduction and Statement of the Theorem}

Let $\Gamma_0(N)$ be the subgroup of elements whose left lower entry is divisible by $N$ in $\mathrm{SL}_2(\mathbb{Z})$ and $\chi$ be a Dirichlet character modulo $N$. Cohen and Oesterl\'{e} \cite{Cohen1976Dimension} stated without proof the dimension formulas for $\mathcal{M}_k(\Gamma_0(N),\chi)$ and $\mathcal{S}_k(\Gamma_0(N),\chi)$, the space of modular forms and that of cusp forms of weight $k$ and character $\chi$ for $\Gamma_0(N)$ respectively. Here $k$ can be integral or half-integral. When $k$ is integral, Quer\cite{Quer} gave a proof of such dimension formulas for all congruence subgroups $\Gamma_H(N)$ sitting between $\Gamma_1(N)$ and $\Gamma_0(N)$. Kaplan and Petrow \cite{Kaplantrace} obtained the $\text{Eichler-Selberg}$ trace formula for a family of congruence
subgroups and then the same dimension formulas of cusp forms for $\Gamma_H(N)$ as proved by Quer \cite{Quer}. Based on such formulas, packages for $\Gamma_0(N)$ by Stein and for $\Gamma_H(N)$ by Quer were implemented in MAGMA and SAGE.

We will be interested on Fuchsian subgroups of $\mathrm{SL}_2(\mathbb{R})$ that normalize $\Gamma_0(N)$. For each $e\parallel N$, that is $e\mid N$ and $\gcd(e,N/e)=1$, there exists $W_e\in\mathrm{SL}_2(\mathbb{R})$ that normalizes $\Gamma_0(N)$ of the form $$W_e=\begin{pmatrix}
a\sqrt{e}&b/\sqrt{e}\\cN/\sqrt{e}&d\sqrt{e}\end{pmatrix},\
a,b,c,d\in \mathbb{Z},\ \det(W_e)=1.$$
Let $\Gamma_0^*(N)$ be the subgroup of $\mathrm{SL}_2(\mathbb{R})$ generated by $\Gamma_0(N)$ and all $W_e$, $e\parallel N$, and $\Gamma_0^+(N)$ be that generated by $\Gamma_0(N)$ and $W_N=\left(\begin{smallmatrix}
0&-1/\sqrt{N}\\\sqrt{N}&0
\end{smallmatrix}\right)$. The resulting groups are independent of the choice of $W_e$ and the quotient group $\Gamma_0^*(N)/\Gamma_0(N)$ is a $2$-elementary group of order $2^{\omega(N)}$, where $\omega(N)$ is the number of distinct prime factors dividing $N$, and $\Gamma_0^+(N)/\Gamma_0(N)\cong \mathbb{Z}/2\mathbb{Z}$. These groups are directly related to Atkin-Lehner theory and see \cite{Atkin-Lehner} for more details. 

Choi and Kim \cite{Choi2013Basis,Choi2019Basis} obtained a canonical basis of the space of weakly holomorphic modular forms with character for the Fricke groups $\Gamma_0^+(p)$ of genus zero by proving the corresponding dimension formulas. As an application, Choi and Im \cite{Choizeros} proved that the zeros in the fundamental domain for $\Gamma_0^+(2)$ of certain weakly holomorphic modular forms lie on the circle with radius $1/\sqrt{2}$. Hanamoto and Kuga \cite{Hanamotozeros} extended this result to $\Gamma_0^+(3)$.

In this paper, we obtain dimension formulas like that of Cohen and Oesterl\'{e}, for the spaces of modular forms of integral weight and with characters for $\Gamma_0^+(N)$ with general $N$.

\vspace{0.3cm}
\noindent
\textbf{Theorem 1.} Let $N>1$, $\chi$ be trivial or a quadratic Dirichlet character modulo $N$, and $\chi^+$ be one of the two characters for $\Gamma_0^+(N)$ induced from $\chi$. For $k\in \mathbb{Z}$ with $\chi(-1)=(-1)^k$,
\begin{align*}
&\text{ dim }\mathcal{S}_k(\Gamma_0^+(N),\chi^+)-\text{ dim }\mathcal{M}_{2-k}(\Gamma_0^+(N),\overline{\chi^+})\\
&=\frac{k-1}{12}\nu_0^+(\chi^+)+\gamma_4(k)\nu_2^+(\chi^+)+
\gamma_3(k)\nu_3^+(\chi^+)+\gamma_8(k)\delta_{N,2}+\gamma_{12}(k)\delta_{N,3}-\frac{\nu_\infty^+(\chi^+)}{2}.
\end{align*}
See the meaning of the notations in Section \ref{section2}. If in addition $N>3$ and $N\equiv 1(\text{mod }4)$ and $\chi$ is quadratic, then Theorem \ref{main1} implies the following $2$-power relation (Corollary \ref{main2}) $$\text{ dim }\mathcal{S}_k(\Gamma_0^+(N),\chi^+)
=\frac{1}{2}\text{ dim }\mathcal{S}_k(\Gamma_0(N),\chi).$$

We will also work on the group $\Gamma^*_0(N)$ in the case when $N$ is square-free and obtain a dimension formula (see Theorem \ref{main3}).
As a corollary, if $N$ is square-free and all prime factors $p$ of $N$ satisfy $p\equiv1(\text{mod }4)$ and $\chi$ is primitive with conductor $N$, then we have the following $2$-power relation
$$\text{ dim }\mathcal{S}_k(\Gamma_0^*(N),\chi^*)
=\frac{1}{2^{\omega(N)}}\text{ dim }\mathcal{S}_k(\Gamma_0(N),\chi).$$
We remark that the $2$-power relation is consistent with the decomposition of $\mathcal{S}_k(\Gamma_0(N),\chi)$ in eigenspaces of the different Atkin-Lehner operators. More precisely, under the isomorphism
\[\mathcal{S}_k(\Gamma_0(N),\chi)=\bigoplus_{\chi^*\text{ lifts }\chi}\mathcal{S}_k(\Gamma_0^*(N),\chi^*),\]
all of the subspaces on the right-hand side have equal dimensions. We also remark that with these dimension formulas one can obtain dimension formulas on Riemann-Roch spaces with prescribed divisor (see (1.6.1) and (3.3.6) of \cite{Shimura1974Onthetrace}). For example, if $k\in \mathbb{Z}, k>2$ and $\Gamma=\Gamma_0(N)$, $\Gamma_0^+(N)$ or $\Gamma_0^*(N)$, then
$$\text{ dim }\mathcal{M}^!_{k,D}(\Gamma,\chi)=\text{ dim }\mathcal{S}_k(\Gamma,\chi)+\text{deg }(D),$$
where $D$ is a cuspidal divisor and $\mathcal{M}^!_{k,D}(\Gamma,\chi)$ is the corresponding Riemann-Roch subspace of weakly holomorphic modular forms (that is, $\text{div } f+D\geq 0$). But we shall not explore such dimensions in this paper.

Shimura \cite{Shimura1974Onthetrace} proved the dimension formula for $\mathcal{S}_k(\Gamma,\chi)$ and $\mathcal{M}_k(\Gamma,\chi)$ for more general Fuchsian groups $\Gamma$ via Riemann-Roch Theorem. From Shimura's formula, we compute explicitly the relevant quantities and obtain above dimension formulas. Such calculation is technical and involves evaluating character sums on elliptic elements. The most tricky part lies on the relevant character sum on the set of extra elliptic points of order $2$, and we employ its connection with integral quadratic forms to work out the desired formula. We then work similarly on $\Gamma_0^*(N)$.

Here is the layout of this paper. In Section \ref{section2}, we first recall Shimura's result and the signature of $\Gamma^+_0(N)$ to prove the dimension formula. Then we explain the connection between elliptic elements in $W_N\Gamma_0(N)$ and quadratic forms of discriminant $-4N$, and by considering generic characters, we prove the corollary on the first $2$-power relation above. In Section \ref{section 3}, we extend above treatment to $\Gamma_0^*(N)$, prove the corresponding dimension formula, and end this paper with a concrete example.

\section{Dimension formula of modular forms with character for $\Gamma_0^+(N)$}\label{section2}

Let $\mathcal{H}$ be the complex upper half plane and put $\mathcal{H}^*=\mathcal{H}\cup \mathbb{Q}\cup \{\infty\}$. For a Fuchsian group $\Gamma$ of the first kind, let $X(\Gamma)=\Gamma\backslash \mathcal{H}^*$ be the corresponding compact Riemann surface and $\pi:\mathcal{H}^*\rightarrow X(\Gamma)$ be the projection. In this paper, $\Gamma$ will be a subgroup of $\mathrm{SL}_2(\mathbb{R})$ generated by $\Gamma_0(N)$ and some matrices $W_e$, $e\parallel N$, of the following form
$$W_e=\begin{pmatrix}
a\sqrt{e}&b/\sqrt{e}\\cN/\sqrt{e}&d\sqrt{e}\end{pmatrix},\
a,b,c,d\in \mathbb{Z},\ \det(W_e)=1.$$
For each $e\parallel N$, the existence of $a,b,c,d$ is clear and their choice makes no difference on the group $\Gamma$, so we may and will require that $b=d=1$. Each $W_e$ normalizes $\Gamma_0(N)$, that is $W_e^{-1}\Gamma_0(N)W_e=\Gamma_0(N)$, hence acts on the spaces of modular forms for $\Gamma_0(N)$, and is called an \emph{Atkin-Lehner involution} (see \cite{Atkin-Lehner}). We denote by $\Gamma_0^e(N)$ the group generated by $\Gamma_0(N)$ and a single $W_e$, and $\Gamma_0^*(N)$ the group generated by $\Gamma_0(N)$ and all possible $W_e$, $e\parallel N$. When $e=N$, set $\Gamma^+_0(N)=\Gamma^N_0(N)$.

For a Dirichlet character modulo $N$, it extends to a character $\chi$ of $\Gamma_0(N)$ via $\chi(\gamma)=\chi(d)$ for $\gamma=\left(\begin{smallmatrix}
a&b\\c&d
\end{smallmatrix}\right)\in \Gamma_0(N)$,
and let $\chi^+$ be one character of $\Gamma_0^+(N)$ satisfying
$\chi^+|_{\Gamma_0(N)}=\chi.$
Therefore $\chi^+(W_N)^2=\chi^+(-I)=\chi(-1)$, and then $\chi^+(W_N)=\pm 1$ if $\chi(-1)=1$ and $\chi^+(W_N)=\pm i$ if  $\chi(-1)=-1$. It can be seen easily that $\chi$ extends to $\Gamma_0^+(N)$ if and only if $\chi$ is quadratic (or trivial), in which case there are exactly two extensions $\chi^+$. Similarly, $\chi$ extends to a character $\chi^*$ of $\Gamma_0^*(N)$ if and only if $\chi$ is quadratic (or trivial), in which case there are exactly $2^{\omega(N)}$ extensions $\chi^*$. Here $\omega(N)$ is equal to the number of distinct prime divisors of $N$.

In this section, we treat the group $\Gamma_0^+(N)$.

\subsection{Dimension formula for $\Gamma_0^+(N)$}

To state the dimension formula for the group $\Gamma_0^+(N)$, we denote by $\nu^+_{0}(N)=\frac{N}{2}\prod\limits_{p\mid N}\left(1+\frac{1}{p}\right)$, $\nu^+_{e}(N)$ and $\nu_{\infty}^+(N)$ the number of order-$e$ elliptic points with $e=2,3,4,6$ and that of cusps for $\Gamma_0^+(N)$ respectively. We introduce 
$$\gamma_3(k)= \left\{
\begin{array}{lcl}
1/3 & & {\text{if } k\equiv0\text{ (mod }3)}\\
0 & & {\text{if } k\equiv1\text{ (mod }3)}\\
-1/3 & & {\text{if } k\equiv2\text{ (mod }3)}\\
\end{array} \right.,\quad  \gamma_4(k)= \left\{
\begin{array}{lcl}
1/4 & & {\text{if } k\equiv0\text{ (mod }4)}\\
-1/4 & & {\text{if } k\equiv2\text{ (mod }4)}\\
i/4 & & {\text{if } k\equiv1\text{ (mod }4)}\\
-i/4 & & {\text{if } k\equiv3\text{ (mod }4)}\\
\end{array} \right.,$$
and  $\gamma_{8}(k)$,  $\gamma_{12}(k)$ in Table \ref{gamma_8(k)} and Table \ref{gamma_12(k)} respectively.
\begin{table}[H]
\centering
\caption{The value of $\gamma_8(k)$}\label{gamma_8(k)}
\vspace{0.3cm}
\begin{tabular}{|c|c|c|c|c|c|c|c|c|}
  \hline
  $k$ (mod 8) &\multicolumn{2}{c|}{0}&\multicolumn{2}{c|}{2}&\multicolumn{2}{c|}{4}
  &\multicolumn{2}{c|}{6}  \\\hline
  $\chi^+(W_2)$ & $+1$ &$-1$ &$+1$ &$-1$ &$+1$ &$-1$ &$+1$ &$-1$ \\\hline
  $\gamma_8(k)$& $\frac{3}{8}$& $-\frac{1}{8}$& $-\frac{3}{8}$& $\frac{1}{8}$& $-\frac{1}{8}$& $\frac{3}{8}$& $\frac{1}{8}$&$-\frac{3}{8}$\\\hline
\end{tabular}
\end{table}

\begin{table}[H]
\centering
\caption{The value of $\gamma_{12}(k)$}\label{gamma_12(k)}
\vspace{0.3cm}
\resizebox{13.5cm}{!}{
\begin{tabular}
{|c|c|c|c|c|c|c|c|c|c|c|c|c|}
\hline \diagbox{$\chi^+(W_3)$}{$k$ (mod $12$)} & $0$ & $1$ & $2$ & $3$ & $4$ & $5$ & $6$ & $7$ & $8$ & $9$ & $10$ & $11$ \\
\hline $+1$ & $\frac{5}{12}$ &  & $-\frac{5}{12}$ &  & $-\frac{3}{12}$ &  & $-\frac{1}{12}$ &  & $\frac{1}{12}$ &  & $\frac{3}{12}$ &  \\
\hline $-1$ & $-\frac{1}{12}$ &  & $\frac{1}{12}$ &  & $\frac{3}{12}$ &  & $\frac{5}{12}$ &  & $-\frac{5}{12}$ &  & $-\frac{3}{12}$ &  \\
\hline $+i$ &  & $-\frac{3}{12}$ &  & $-\frac{1}{12}$ &  & $\frac{1}{12}$ & & $\frac{3}{12}$ &  & $\frac{5}{12}$ &  & $-\frac{5}{12}$ \\
\hline $-i$ &  & $\frac{3}{12}$ &  & $\frac{5}{12}$ &  & $-\frac{5}{12}$ &  & $-\frac{3}{12}$ &  & $-\frac{1}{12}$ &  & $\frac{1}{12}$ \\
\hline
\end{tabular}}
\end{table}

We remark that $\gamma_4(k)$, $\gamma_3(k)$, $\gamma_8(k)$ and $\gamma_{12}(k)$ correspond to elliptic points of order $2$, $3$, $4$ and $6$ respectively. Order-$4$ and order-$6$ elliptic points occur precisely when $N = 2$ and $N = 3$ respectively. Since $\chi$ is trivial when $N=2$, $k$ must be even and odd $k$ does not appear in Table \ref{gamma_8(k)}, while $\chi$ can be trivial or quadratic when $N=3$ and both parities of $k$ appear in Table \ref{gamma_12(k)}.
\begin{Thm}\label{main1}
Let $N>1$ be an integer, and $\chi$ be trivial or a quadratic Dirichlet character modulo $N$ that extends to a character $\chi^+$ of $\Gamma_0^+(N)$. For $k\in \mathbb{Z}$ with $\chi(-1)=(-1)^k$,
\begin{align}
&\text{ dim }\mathcal{S}_k(\Gamma_0^+(N),\chi^+)-\text{ dim }\mathcal{M}_{2-k}(\Gamma_0^+(N),\overline{\chi^+})\label{Gamma0N+equ}\\
&=\frac{k-1}{12}\nu_0^+(N)+\gamma_4(k)\nu_2^+(\chi^+)+
\gamma_3(k)\nu_3^+(\chi^+)+\gamma_8(k)\delta_{N,2}+\gamma_{12}(k)\delta_{N,3}-\frac{\nu_\infty^+(\chi^+)}{2},\notag
\end{align}
where $\delta_{a,b}=1$ if $a=b$ and $0$ otherwise,
\[\nu_2^+(\chi^+)=\sum_{z\in E^+_{2}}\chi^+(\gamma_z),\quad \nu_3^+(\chi^+)=\sum_{s\in A_3'(N)}\chi(s),\]
\[\nu_\infty^+(\chi^{+})=\left\{
\begin{array}{lcl}
\frac{1}{2}\sum\limits_{\substack{c\mid N \\ \gcd(c,N/c)\mid N/f}}\varphi\left(\gcd\left(c,N/c\right)\right) & & {N>1,\ N\neq4}\\
2 & & {N=4,\ \chi^+ \text{ is trivial}}\\
1 & & {N=4,\ \chi^+(W_N)=-1}\\
3/2 & & {N=4,\ \chi^+(W_N)=i}\\
1/2 & & {N=4,\ \chi^+(W_N)=-i}\\
\end{array}\right.,\]
$E^+_{e}$ is the set of elliptic points of order $e$ for $\Gamma_0^+(N)$, $A_3'(N)=\{s\in\Z/N\Z|s^2+s+1\equiv 0(\text{mod }N)\}/\sim$
with $s\sim s'$ if $s+s'\equiv-1 (\text{mod }N)$, $\gamma_z$ is one generator of $\Gamma^+_0(N)_{\tau}/\{\pm I\}$ with $z=\pi(\tau)\in X(\Gamma_0^+(N))$,   $\varphi$ is Euler's totient function, and $f$ is the conductor of $\chi$.
\end{Thm}
%\subsection{Proof of Theorem \ref{main1}}
\begin{proof} Set $\zeta_m=e^{2\pi i/m}$.
Applying formula (3.3.6) in \cite{Shimura1974Onthetrace}, we have
 \begin{align}
  &\text{ dim }\mathcal{S}_k(\Gamma_0^+(N),\chi^+)-\text{ dim }\mathcal{M}_{2-k}(\Gamma_0^+(N),\overline{\chi^+})\label{Shimuragamma0N+}\\
&=\frac{k-1}{12}\nu_0^+(N)-\frac{\nu_2^+(N)}{4}-\frac{\nu_3^+(N)}{3}-\frac{3\nu_4^+(N)}{8}-\frac{5\nu_6^+(N)}{12}
-\frac{\nu_{\infty}^+(N)}{2}+\sum_{z\in \mathfrak{R}^+_{N}}\mu_z',\notag
\end{align}
where $\mathfrak{R}^+_{N}$ consists of all elliptic points and cusps for $X(\Gamma_0^+(N))$ and $\mu_z'\in\mathbb Q$ satisfies
$$\mu_z'\equiv v_z(g) \text{ mod }\mathbb Z,\quad 0\leq\mu_z'<1 \text{ for } 0\neq g\in \mathcal{M}_{2-k}(\Gamma_0^+(N),\overline{\chi^+})$$ with $v_z(g)=n/|\Gamma^+_0(N)_{\tau}/\{\pm I\}|$, where $n$ is the order of Fourier expansion of $g$ at $\tau$. See \cite{Shimura1974Onthetrace} for more details.

%\noindent
Recall that when $N=2$, $\Gamma^{+}_0(2)$ has two elliptic points,  $\frac{i}{\sqrt{2}}$ of order $2$ and $\frac{1}{\sqrt{2}}\zeta_8^3$ of order $4$ (see \cite{Tsuyoshiellipticpoint}). The elliptic points of $\Gamma^{+}_0(3)$ are $\frac{i}{\sqrt{3}}$ of order $2$ and $\frac{1}{\sqrt{3}}\zeta_{12}^5$ of order $6$ (see also \cite{Tsuyoshiellipticpoint}).  When $N\geq4$, any elliptic point of $\Gamma^+_0(N)$ is of order $2$  or $3$ (see \cite{Fricke}, p.357-367, for more details). Therefore, by comparing the equations \eqref{Gamma0N+equ} and \eqref{Shimuragamma0N+}, we need to prove that when $N>3$,
$$\sum\limits_{z\in E^+_{2}}\left(\mu_z'-\frac{1}{4}\right)=\gamma_4(k)\nu^+_2(\chi^+), \sum\limits_{z\in E^+_{3}}\left(\mu_z'-\frac{1}{3}\right)=\gamma_3(k)\nu^+_3(\chi^+), \sum\limits_{\frac{a}{c}\in P^+_{N}}\mu'_{\frac{a}{c}}=\frac{\nu_{\infty}^+(N)}{2}-\frac{\nu^+_{\infty}
(\chi^+)}{2}$$
with $P^+_{N}$ being the set of cusps for $\Gamma_0^+(N)$. If $N=2$ or $3$, we have to include order-$4$ and order-$6$ elliptic points and prove that
$$\sum\limits_{z\in{E^+_{4}}}\left(\mu_z'-\frac{3}{8}\right)=\gamma_8(k)\delta_{N,2},\quad \sum\limits_{z\in{E^+_{6}}}\left(\mu_z'-\frac{5}{12}\right)=\gamma_{12}(k)\delta_{N,3}. $$
In the following we verify these formulas on elliptic points first and then that on cusps.

\medskip
\paragraph*{\bf (1) Computation on elliptic points.}
We begin with the computation on order-$2$ elliptic points. Let $z=\pi(\tau)$ be an elliptic point of order $2$ and then $\mu_z'\in [0,1)$ is the unique rational number such that $e^{2\pi i \mu'_z}=i^{2-k}{\chi^+}(\gamma_z)^{-1}$. When $k$ is even, $\chi(-1)=1$, so $\chi^+$ is a real character with $\chi^+(W_N)=\pm1$ and $e^{2\pi i \mu'_z}=i^{2-k}\chi^+(\gamma_z)$. If $k\equiv 0(\text{mod }4)$, then $e^{2\pi i \mu'_z}=-\chi^+(\gamma_z)$ and  $\chi^+(\gamma_z)=\pm 1$ if and only if $\mu'_z=\frac{1}{2},0$ respectively, so $\mu'_z-\frac{1}{4}=\frac{1}{4}\chi^+(\gamma_z)$. Similarly, we can get $\mu'_z-\frac{1}{4}=-\frac{1}{4}\chi^+(\gamma_z)$ when $k\equiv 2(\text{mod }4)$. On the other hand, when $k$ is odd, $\chi(-1)=-1$ and $\chi^+(W_N)=\pm i$. We claim that either $4\mid N$ or there exists a prime divisor $p\equiv 3(\text{mod 4})$ of $N$, in which case it is well-known that $\Gamma_0(N)$ has no elliptic points of order $2$.  Indeed, assume that $4\nmid N$. The $2$-component $\chi_2$ of $\chi$ is trivial and hence $\chi_2(-1)=1$. Moreover for each prime divisor $p\equiv 1(\text{mod }4)$ of $N$, the quadratic $p$-component $\chi_p$ must satisfy $\chi_p(-1)=1$, and the claim follows from that $\chi(-1)=-1$. Therefore we just need to compute $\sum_{z\in E^+_{2}}(\mu_z'-\frac{1}{4}),$ where $z$ is fixed by certain elements in $W_N\Gamma_0(N)$. If $k\equiv1(\text{mod }4),\ e^{2\pi i \mu'_z}=i{\chi^+}(\gamma_z)^{-1}$, so $\mu_z'=0,\frac{1}{2}$ if and only if $\chi^+(\gamma_z)=\pm i$ respectively. Hence $\sum_{z\in E^+_{2}}(\mu_z'-\frac{1}{4})=\frac{i}{4}\nu_2^+(\chi^+)$. The case $k\equiv3(\text{mod }4)$ is similar and we get that $\sum_{z\in E^+_{2}}(\mu_z'-\frac{1}{4})=-\frac{i}{4}\nu_2^+(\chi^+)$. Putting everything together, we obtain the formula  $\sum_{z\in E^+_{2}}(\mu_z'-\frac{1}{4})=\gamma_4(k)\nu^+_2(\chi^+)$ as desired.

Let $z=\pi(\tau)$ be an order-$3$ elliptic point and $\gamma_z$ be a generator of $\Gamma^+_0(N)_{\tau}/\{\pm I\}$. Clearly, $\gamma_z^2$ is also a generator and $\gamma_z^2\in\Gamma_0(N)$, so $\gamma_z\in\Gamma_0(N)$ and $z$ is an elliptic point of $\Gamma_0(N)$. Therefore, we can choose $\tau=\frac{s+\zeta_6}{s^2+s+1}$ and $\gamma_z=\left(\begin{smallmatrix} s & -1 \\ s^2+s+1 & -s-1\end{smallmatrix}\right)$ (see \cite{Diamond2005}, p.96-97), where $s$ runs through all solutions to $s^2+s+1\equiv 0 (\text{mod } N)$ with $s$ and $-s-1$ paired into one orbit via $W_N$. In particular, $E^+_{3}$ bijects to $A'_3(N)$ and the summation over $E^+_{3}$ equals that over $A'_3(N)$. Now $e^{2\pi i \mu'_z}=\zeta_3^{2-k}\overline{\chi^+(\gamma_z)}=\zeta_3^{2-k}\overline{\chi(s^2)}=\zeta_3^{2-k}$ since $\chi$ is quadratic, and by elementary computation we obtain $\sum_{z\in E^+_{3}}(\mu_z'-\frac{1}{3})=\gamma_3(k)\nu^+_3(\chi^+)$.

Finally we treat the order-$4$ and order-$6$ elliptic points and as we have already seen they exist precisely when $N=2$ and $N=3$ respectively. Since $\Gamma^+_0(2)$ has a unique order-$4$ elliptic point $z=\pi(-\frac{1}{2}+\frac{i}{2})$, we may choose $\gamma_{z}$ to be $\left(\begin{smallmatrix} -\sqrt{2} & -1/\sqrt{2} \\ \sqrt{2} & 0\end{smallmatrix}\right)$, so $e^{2\pi i\mu_{z}'}=\zeta_8^{6-3k}\chi^+(\gamma_z)^{-1}$. Similarly, $\Gamma^+_0(3)$ has a unique order-$6$ elliptic point $z=\pi(-\frac{1}{2}+\frac{\sqrt{3}i}{6})$ and then $e^{2\pi i\mu_{z}'}=\zeta_{12}^{10-5k}\chi^+(\gamma_z)^{-1}$ with $\gamma_{z}=\left(\begin{smallmatrix} -\sqrt{3} & -1/\sqrt{3} \\ \sqrt{3} & 0\end{smallmatrix}\right)$. It follows that $\sum_{z\in{E^+_{4}}}\left(\mu_z'-\frac{3}{8}\right)=\gamma_8(k)$ when $N=2$ and  $\sum_{z\in{E^+_{6}}}\left(\mu_z'-\frac{5}{12}\right)=\gamma_{12}(k)$ when $N=3$.

\medskip
\paragraph*{\bf (2) Computation on cusps.}

Recall that a set of representatives of  $\Gamma_0(N)$-orbits of cusps is
$$\left\{\frac{a}{c}\colon\ c|N, \ a \imod \gcd(c,N/c), \ \gcd(a,c,N/c)=1\right\}$$
and no extra cusps for $\Gamma_0^+(N)$ appear since $\Gamma^+_0(N)$ is commensurable with $\Gamma_0(N)$.

Assume first that $N\neq4$. In this case, $W_N$ pairs distinct cusps of $\Gamma_0(N)$, that is $W_N\frac{a}{c}=\frac{a'}{N/c}$, where $a'\equiv -a (\text{mod }\gcd(c,N/c))$. For a cusp $a/c$ of $\Gamma^+_0(N)$, the corresponding parabolic element $\gamma_{a/c}$ can be chosen to be $\left(\begin{smallmatrix} 1-ach & a^2h \\ -c^2h & 1+ach\end{smallmatrix}\right)\in \Gamma_0(N)$ with $h=N/\gcd(N,c^2)$. So $e^{2\pi i\mu_{a/c}'}={\chi^+}(\gamma_{a/c})^{-1}=\chi(1+ach)$ and it suffices to prove that for any $c\mid N, c\leq\sqrt{N},$
$$ \sum\limits_{a\in \left(\mathbb{Z}/\gcd(c,N/c)\mathbb{Z}\right)^{\times}}\mu'_{a/c}=\left\{
\begin{array}{lcl}
0 & & {\text{if } f\mid \frac{N}{\gcd(c,N/c)}}\\
\frac{1}{2}\phi(\gcd(c,N/c)) & & {\text{if } f\nmid\frac{N}{\gcd(c,N/c)}}\\
\end{array} \right., $$
where $f$ is the conductor of $\chi$.

If $f\mid\frac{N}{\gcd(c,N/c)}$, then $ach\equiv0\text{(mod }f)$ and $\mu_{a/c}'=0$ for any $a\in \left(\mathbb{Z}/\gcd(c,N/c)\mathbb{Z}\right)^{\times}$. If $f\nmid\frac{N}{\gcd(c,N/c)}$, we claim that $\chi(1+ach)\neq 1$ for any $a\in\left(\mathbb{Z}/\gcd(c,N/c)\mathbb{Z}\right)^{\times}$. Indeed, if $\chi(1+ach)=1$ for some $a$, then \[\chi\left(\left(1+a\frac{N}{\gcd\left(c,N/c\right)}\right)^n\right)=\chi\left(1+\frac{N}{\gcd\left(c,N/c\right)}\right)=1\] for any $n\in \mathbb{N}$ such that $an \equiv1(\text{mod }\gcd\left(c,N/c\right))$ since $\left(1,\gcd\left(c,N/c\right)\right)=1$, and it follows that
$\chi\left(1+\frac{N}{\gcd\left(c,N/c\right)}\mathbb{Z}\right)=1$.
Therefore $f$ divides $\frac{N}{\gcd\left(c,N/c\right)}$, which contradicts to the condition $f\nmid\frac{N}{\gcd(c,N/c)}.$ Clearly, $\gcd\left(c,N/c\right)>1,$ since otherwise $\chi(1+ach)=1.$ If $\gcd\left(c,N/c\right)=2,$ then there is only one cusp $1/c$ in this case and $\chi(1+ach)^{2}=\chi(1+2ach)=1.$ It follows that $\chi(1+ach)=-1$ and $\mu_{1/c}'=\frac{1}{2}$, which verifies the formula. If $\gcd\left(c,N/c\right)>2,$ then $a/c$ and $-a/c$ are distinct cusps and $\chi(1+ach)\cdot \chi(1-ach)=1$. So $\mu_{a/c}'+\mu_{-a/c}'=1$, since $\chi(1+ach)=e^{2\pi i \mu'_{a/c}}\neq 1$ and hence $0<\mu'_{a/c}, \mu'_{-a/c}<1$. Therefore
$$ \sum\limits_{a\in \left(\mathbb{Z}/\gcd(c,N/c)\mathbb{Z}\right)^{\times}}\mu'_{a/c}=\frac{1}{2}\sum\limits_{a\in \left(\mathbb{Z}/\gcd(c,N/c)\mathbb{Z}\right)^{\times}}
(\mu_{a/c}'+\mu_{-a/c}')=\frac{1}{2}\phi(\gcd(c,N/c)).$$

Now we treat the case when $N=4$. Recall that $W_4$ acts on the set of cusps of $\Gamma_0(4)$ and it switches $0$ and $\infty$ but fixes $1/2$. For cusps $1/2$ and $\infty$, we get $e^{2\pi i \mu'_{1/2}}=\chi^+(\gamma_{1/2})^{-1}$ with $\gamma_{1/2}=\left(\begin{smallmatrix}0&1/2\\-2&2\end{smallmatrix}\right)$ and $e^{2\pi i \mu'_{\infty}}=\chi(1+ach)$ respectively. Elementary computation shows that $\sum_{a/c\in P^+_{N}}\mu'_{a/c}=\frac{1}{2}\nu_{\infty}^+(N)-\frac{1}{2}\nu^+_{\infty}
(\chi^+)$. This completes the proof.
\end{proof}

We remark that similar but easier computation as in the proof of Theorem \ref{main1} gives a proof of Cohen and Oesterl\'{e}'s formula in \cite{Cohen1976Dimension}.  In Theorem \ref{main1}, the only term to be determined explicitly is $\nu_2^+(\chi^+)$, which will be done soon in Proposition \ref{helpproofcor1} below. Before that, we record the dimension formulas for the subspace of cusp forms and that for the Eisenstein subspace.

\begin{Cor}
For $k=2$ and $\chi^+$ is non-trivial or $k>2$, we have
\begin{align*}
\text{ dim }\mathcal{S}_k(\Gamma^+_0(N),\chi^+)&=
\frac{k-1}{12}\nu_0^+(N)+\gamma_4(k)\nu_2^+(\chi^+)+
\gamma_3(k)\nu_3^+(\chi^+)\\
&+\gamma_8(k)\delta_{N,2}+\gamma_{12}(k)\delta_{N,3}-\frac{\nu_\infty^+(\chi^+)}{2}. \end{align*}
If $k=2$ and $\chi^+$ is trivial, then
$\text{dim }\mathcal{S}_2(\Gamma^+_0(N),\chi^+)=g(\Gamma_0^+(N))$ with the genus
\[g(\Gamma_0^+(N))=\frac{\nu_0^+(N)}{12}-\frac{\nu_2^+
(N)}{4}-\frac{\nu_3^+(N)}{3}-\frac{3}{8}\delta_{N,2}-\frac{5}{12}\delta_{N,3}-\frac{\nu_\infty^+(N)}{2}+1.\]
\end{Cor}
\begin{proof} When $k\geq2$, it is obvious that
\begin{align}\label{dimen}
\text{ dim }\mathcal{M}_{2-k}(\Gamma_0^+(N),\overline{\chi^+})= \left\{
\begin{array}{lcl}
1 & & {\text{if } k=2 \text{ and }\chi^+ \text{ is trivial}}\\
0 & & {\text{if } k>2 \text{ or } k=2 \text{ and }\chi^+ \text{ is non-trivial}}\\
\end{array} \right.
\end{align}
and the dimension formula for $\mathcal{S}_k(\Gamma^+_0(N),\chi^+)$ follows by Theorem \ref{main1}.
\end{proof}
\begin{Cor}
Let us denote by $\mathcal{E}_k(\Gamma_0^+(N),\chi^+)$ the Eisenstein subspace, namely the orthogonal complement of $\mathcal{S}_k(\Gamma_0^+(N),\chi^+)$ in $\mathcal{M}_k(\Gamma_0^+(N),\chi^+)$. For even integer $k\geq 2$, we have $$\text{dim }\mathcal{E}_k(\Gamma_0^+(N),\chi^+)=\left\{
\begin{array}{lcl}
\nu_{\infty}^+(\chi^+)-1 & & {\text{if } k=2 \text{ and }\chi^+ \text{ is trivial}}\\
\nu_{\infty}^+(\chi^+) & & {\text{if } k>2 \text{ or } k=2 \text{ and }\chi^+ \text{ is non-trivial}}\\
\end{array} \right..$$
\end{Cor}
\begin{proof}
When $k$ is even, $\chi^+$ is a real character. If $k\geq2$, formula \eqref{Gamma0N+equ} becomes
\begin{align*}
&\text{ dim }\mathcal{S}_k(\Gamma_0^+(N),\chi^+)-\text{ dim }\mathcal{M}_{2-k}(\Gamma_0^+(N),\chi^+)\\&=\frac{k-1}{12}\nu_0^+(N)+\gamma_4(k)\nu_2^+(\chi^+)+
\gamma_3(k)\nu_3^+(\chi^+)+\gamma_8(k)\delta_{N,2}+\gamma_{12}(k)\delta_{N,3}-\frac{\nu_\infty^+(\chi^+)}{2}
\end{align*} and by interchanging $k$ and $2-k$, we have
\begin{align*}
\text{ dim }\mathcal{M}_k(\Gamma_0^+(N),\chi^+)&=\frac{k-1}{12}\nu_0^+(N)+\gamma_4(k)\nu_2^+(\chi^+)\\&+
\gamma_3(k)\nu_3^+(\chi^+)+\gamma_8(k)\delta_{N,2}+\gamma_{12}(k)\delta_{N,3}+\frac{\nu_\infty^+(\chi^+)}{2},
\end{align*} since $\gamma_e(k)=-\gamma_e(2-k)$.
By subtracting the two equations, the corollary then follows from formula \eqref{dimen}.
\end{proof}

\subsection{Elliptic elements for $\Gamma_0^+(N)$ and quadratic forms}

In the dimension formula of Theorem \ref{main1}, the term $\nu_2^+(\chi^+)$ is not fully explicit since the set $E^+_{2}$ of elliptic points of order $2$ has not been described. In this subsection, we employ integral quadratic forms to characterize $E^+_{2}$ and obtain the formula of $\nu_2^+(\chi^+)$. As an application, when $N\equiv 1 (\text{mod }4)$, we prove a $2$-power relation on dimensions.

We first define $E^+_{2}=E^{+}_{2,0}\bigsqcup E^{+}_{2,1}$, where $E^{+}_{2,0}$ and $E^{+}_{2,1}$ consist of order-$2$ elliptic points fixed by elements of $\Gamma_0(N)$ and $W_N\Gamma_0(N)$ respectively.
Let $G^{+}_{0}$ and $G^{+}_{1}$ be the set of generators of $\Gamma^+_0(N)_\tau/\{\pm I\}$, one for each $z=\pi(\tau)\in E^{+}_{2,0}$ and $E^{+}_{2,1}$ respectively, and we define $G^+=G^{+}_{0}\bigsqcup G^{+}_{1}$.

We next recall generic characters of discriminant $-4N$, and then characterize $E^{+}_{2}$ by the connection between elliptic points fixed by $W_N\Gamma_0(N)$ and quadratic forms of discriminant $-4N$ or $-N$ given by Fricke\cite{Fricke}. For more details on generic characters, one may consult Chapter $4$ of \cite{binaryform}.

Define $\mathscr{Q}_{\Delta}^0$ the set of all primitive quadratic forms $(a,b,c)$ of discriminant $\Delta=b^2-4ac$, and $h(\Delta)$ the number of $\mathrm{SL}_2(\Z)$-equivalence classes of primitive, positive definite binary quadratic forms of discriminant $\Delta$. Let $\Delta'$ be the corresponding fundamental discriminant, namely $\Delta'\mid\Delta$ and $\Delta/\Delta'\equiv 0$ or $1(\text{mod }  4)$. The genus character associated to $(\Delta,\Delta')$ is a $\mathrm{SL}_2(\Z)$-invariant function $\chi_{\Delta'}: \mathscr{Q}_{\Delta}^0 \rightarrow\{\pm 1\}$ such that  $\chi_{\Delta'}((a,b,c))=\left(\frac{\Delta'}{r}\right)$ for any integer $r$ prime to $\Delta'$ and represented by $(a,b,c)$. Such $r$ always exists and the value of $\left(\frac{\Delta'}{r}\right)$ is independent of its choice. We can decompose $\left(\frac{\Delta}{\cdot}\right)$ to obtain the so-called generic characters; explicitly, let $v_p(N)$ be the $p$-adic valuation of $N$ and the generic characters of the discriminant $\Delta=-4N$ are given in Table \ref{generic character} (see Page $52$ in \cite{binaryform}). Observe that the genus character $\chi_{\Delta'}$ is exactly the product of all generic characters of $\Delta'$.

\begin{table}
\centering
\caption {Generic characters of $\Delta=-4N$}\label{generic character}
%\vspace{0.3cm}
\begin{tabular}{|c|c|c|c|}
\hline
  \multicolumn{2}{|c|}{$\text{$N$(mod 4)}$} & \multicolumn{2}{|c|}{$\text{Generic characters of $-4N$}$}   \\\hline
  \multirow{2}*{$0$} & $v_2(N)=2$ & \multirow{6}*{$\quad\left(\frac{\cdot}{p}\right),\ p\mid N\quad$} &  $\left(\frac{-1}{\cdot}\right)$ \\\cline{2-2} \cline{4-4}
  & $v_2(N)\geq3$ & &  $\left(\frac{-1}{\cdot}\right)\ \left(\frac{2}{\cdot}\right) $ \\\cline{1-2} \cline{4-4}
  \multicolumn{2}{|c|}{$1$} & & $\left(\frac{-1}{\cdot}\right)$ \\\cline{1-2} \cline{4-4}
  \multirow{2}*{$\qquad2\qquad$} & $2(\text{mod }8)$ & & $\left(\frac{-2}{\cdot}\right)$ \\ \cline{2-2} \cline{4-4}
  & $6(\text{mod }8)$ & & $\left(\frac{2}{\cdot}\right)$ \\\cline{1-2} \cline{4-4}
  \multicolumn{2}{|c|}{$3$} & &  \\\hline
\end{tabular}
\end{table}

Let $r$ be odd such that $\gcd(r,\Delta)=1$ and represented by a quadratic form of discriminant $\Delta$. Choose a fixed order on the set of these generic characters and they can be considered as a vector-valued function from integer $r$ to an tuple with entries either $+1$ or $-1$.
Recall that two integral quadratic forms are in the same genus if they are equivalent over $\R$ and over $\mathbb Z_p$ for every prime $p$. For binary quadratic forms, a genus is the set of classes of quadratic forms that possesses the same assigned values of generic characters, and the assigned values of generic characters are independent of the choice of $r$ represented by a quadratic form in a given genus (see Proposition 4.3 of\cite{binaryform}).

When $N\equiv 0,1,2(\text{mod }4)$, Fricke\cite{Fricke} showed that there is a $1$-$1$ correspondence between $G^{+}_{1}$ and $\mathscr{Q}_{-4N}^0/\mathrm{SL}_2(\Z)$ given by $\gamma=\left(
   \begin{smallmatrix}
     a\sqrt{N} & -b/\sqrt{N} \\
     c\sqrt{N} & -a\sqrt{N} \\
   \end{smallmatrix}
 \right)\mapsto(cN,-2aN,b)$. Note that either both $b$ and $c$ are odd, or one is odd and another one is even, so we may assume that $\gamma$ satisfies $\gcd(b,-4N)=1$ by changing the representatives of elliptic elements of $\Gamma^+_0(N)$. Thus we can take $r=b$ and decompose $G^{+}_{1}$ by the assigned values of generic characters of $\Delta=-4N$ at $b$ (see Example \ref{classify}).

When $N\equiv 3 (\text{mod }4)$, Fricke showed also two $1$-$1$ correspondences, between ${(G^{+}_{1})}^{\prime}$ and $\mathscr{Q}_{-4N}^0/\mathrm{SL}_2(\Z)$, and between ${(G^{+}_{1})}^{\prime\prime}$ and $\mathscr{Q}_{-N}^0/\mathrm{SL}_2(\Z)$, where ${(G^{+}_{1})}^{\prime}$ is the subset of $G^{+}_{1}$ with element of the form $\left(
   \begin{smallmatrix}
     a\sqrt{N} & -b/\sqrt{N} \\
     c\sqrt{N} & -a\sqrt{N} \\
   \end{smallmatrix}
 \right)$ satisfying $\gcd(cN,-2aN,b)=1$ and $\gcd(b,-4N)=1$, and ${(G^{+}_{1})}^{\prime\prime}$ is the complement of ${(G^{+}_{1})}^{\prime}$ in $G^{+}_{1}$ whose element satisfies $\gcd(cN,-2aN,b)=2$ and $2\parallel b$. Similarly, we take $r=b$ to decompose ${(G^{+}_{1})}^{\prime}$ and $r=b/2$ to decompose ${(G^{+}_{1})}^{\prime\prime}$.

\begin{table}
\centering
\caption{Equivalence classes in each genus when $N=65$}\label{example N=65 decomposition of quadratic}
%\vspace{0.3cm}
\setlength{\abovecaptionskip}{2pt}
\setlength{\belowcaptionskip}{2pt}
\begin{tabular}{|ccc|c|}
\hline
  $\left(\frac{-1}{\cdot}\right)$ & $\left(\frac{\cdot}{5}\right)$ & $\left(\frac{\cdot}{13}\right)$ &
  $(cN,-2aN,b)$ of discriminant $-4N$   \\\hline
  $+$ & $+$ & $+$ &  $(N,0,1)$ $(29N,-4N,9)$ \\\hline
  $+$ & $-$ & $-$ &  $(2N,-2N,33)$ $(18N,-14N,177)$ \\\hline
  $-$ & $+$ & $-$ &  $(6N,-2N,11)$ $(6N,2N,11)$ \\\hline
  $-$ & $-$ & $+$ &  $(22N,2N,3)$ $(22N,-2N,3)$ \\
  \hline
\end{tabular}
\end{table}

\begin{Exa}\label{classify}
When $N=65$, the generic characters of $\Delta=-4N$ and the equivalence classes of each genus are listed in Table \ref{example N=65 decomposition of quadratic}.
\end{Exa}

Recall that $h(\Delta)$ is the number of $\mathrm{SL}_2(\Z)$-equivalence classes of primitive, positive definite binary quadratic forms of discriminant $\Delta$. Based on the above results, we provide the explicit formula of $\nu^+_2(\chi^+)$.

\begin{Prop}\label{helpproofcor1}
For quadratic character $\chi$ modulo $N$ with conductor $f$, we decompose $\chi$ into $p$-components as $\chi=\prod_{p|N}\chi_p$. For $-4N=N'l^2$ with $N\geq5$ and $N'$ being the uniquely determined fundamental discriminant, we have
$$\nu_2^+(\chi^+)=\left\{\begin{array}{lcl}
\frac{\nu_2(\chi)}{2}+h(-4N)\chi^+(W_N) & & {\text{if } v_2(N)\geq2 \text{ is even and }f= -N'}\\
\frac{\nu_2(\chi)}{2} & & {\text{if }v_2(N)\geq2\text{ is odd, }f=-N'\text{ and } \chi_2\neq\left(\frac{-2}{\cdot}\right)}\\
\frac{\nu_2(\chi)}{2}+h(-4N)\chi^+(W_N) & & {\text{if } v_2(N)\geq2 \text{ is odd, } f=-N' \text{ and }\chi_2=\left(\frac{-2}{\cdot}\right)}\\
\frac{\nu_2(\chi)}{2}+h^-\chi^+(W_N) & & {\text{if } N\equiv 3(\text{mod }4),f=-N'\text{ with } f\equiv 3(\text{mod }8)}\\
\frac{\nu_2(\chi)}{2}+h^+\chi^+(W_N) & & {\text{if } N\equiv 3(\text{mod }4),f=-N'\text{ with }f\equiv 7(\text{mod }8) }\\
\frac{\nu_2(\chi)}{2} & & {\text{otherwise }}\\
\end{array} \right.,$$
where $\nu_2(\chi)=\sum_{s\in A_4(N)}\chi(s)$, $A_4(N)=\{s\in\mathbb{Z}/N\mathbb{Z}:s^2+1\equiv 0(\text{mod } N)\}$,  $h^+=h(-4N)+h(-N)$, and $h^-=h(-4N)-h(-N)$.
\end{Prop}
\begin{proof}
Write
$$\nu_2^+(\chi^+)=\sum\limits_{\gamma_z\in G^{+}}\chi^+(\gamma_z)=\sum\limits_{\gamma_z\in G_{0}^{+}}\chi^+(\gamma_z)+\sum\limits_{\gamma_z\in G_{1}^{+}}\chi^+(\gamma_z)$$ and we begin with the computation on the first term $\sum_{\gamma_z\in G_{0}^{+}}\chi^+(\gamma_z)$.
Recall that $W_N$ pairs distinct $\Gamma_0(N)$-orbits of order-$2$ elliptic points $\tau$ fixed by $\gamma$ and $W_N\tau$ fixed by $W_N\gamma W_N^{-1}$, then $\chi(\gamma)=\chi(W_N\gamma W_N^{-1})$ and $$\sum_{\gamma_z\in G_{0}^{+}}\chi^+(\gamma_z)=\sum_{\gamma_z\in G} \frac{\chi(\gamma_z)}{2}=\frac{\nu_2(\chi)}{2},$$ where
$G$ is the set of $\Gamma_0(N)$-orbits of elliptic elements of order $2$ and bijects to $A_4(N)$.
Now we treat the second term $\sum_{\gamma_z\in G_{1}^{+}}\chi^+(\gamma_z)$ case by case and the statement follows accordingly.
\medskip
\paragraph*{\bf (1) $N\equiv 1(\text{mod }4)$.}
 We first assume that $-4N$ is a fundamental discriminant, in which case the generic characters are $\left(\frac{\cdot}{p}\right)$, $p\mid N$ and $\left(\frac{-1}{\cdot}\right)$. Define
 $$R_{f,+}=\left\{\left(\begin{array}{cc}a\sqrt{N} & -b/\sqrt{N} \\c\sqrt{N} & -a\sqrt{N} \\
 \end{array}\right)
\in G^{+}_{1}\Bigg|\left(\frac{b}{f}\right)=+1\right\}$$
and
$$R_{f,-}=\left\{\left(
\begin{array}{cc}
a\sqrt{N} & -b/\sqrt{N} \\
c\sqrt{N} & -a\sqrt{N} \\
\end{array}
\right)
\in G^{+}_{1}\Bigg|\left(\frac{b}{f}\right)=-1\right\}.$$
Recall that there is a $1$-$1$ correspondence between $G^{+}_{1}$ and $\mathscr{Q}_{-4N}^0/\mathrm{SL}_2(\Z)$. Theorem $4.16$ and Proposition $4.18$ in \cite{binaryform} show that all the genera have the same number of quadratic forms for any discriminant, and in particular if the discriminant is fundamental, the product of assigned values of generic characters for any given genus is $+1$. Hence the cardinality of $R_{f,+}$ is equal to that of $R_{f,-}$ for $f<N$  and the above conclusion still holds for $f=N$ since $\left(\frac{-1}{\cdot}\right)$ is included in the set of generic characters of $-4N$. So
\begin{align*}
\sum_{\gamma_z\in G^{+}_{1}}\chi^+(\gamma_z)=\chi^+(W_N)\sum_{\gamma_z\in G^{+}_{1} }\chi(b)=\chi^+(W_N)\sum_{\gamma_z\in R_{f,+}\bigsqcup R_{f,-}}\left(\frac{b}{f}\right)=0
\end{align*}
with $\gamma_z=\left(\begin{smallmatrix}
 a\sqrt{N} & -b/\sqrt{N} \\c\sqrt{N} & -a\sqrt{N} \\\end{smallmatrix}\right)=W_N\left(\begin{smallmatrix}
 c & -a \\-aN & b \\\end{smallmatrix}\right)$, $\gcd(b,-4N)=1$.

If $-4N$ is not a fundamental discriminant, write $-4N=N'l^2$ with $4\parallel N'$ and $2\nmid l$. If all prime factors $p$ of $l$ also divide $N'$, the set of generic characters of discriminant $-4N$ equals that of fundamental discriminant $N'$. Applying Theorem $4.16$ and Proposition $4.18$ in \cite{binaryform} again, we have $\sum_{\gamma_z\in G^{+}_{1}}\chi^+(\gamma_z)=0$. Otherwise, $\left(\frac{\cdot}{p}\right)$ is added in the set of generic characters of discriminant $-4N$, doubling the number of genera of discriminant $N'$, and we still have $\sum_{\gamma_z\in G^{+}_{1}}\chi^+(\gamma_z)=0$.
\medskip
\paragraph*{\bf (2) $N\equiv 2(\text{mod }4)$.}
Write $-4N=N'l^2$ with $8\parallel N'$ and $2\nmid l$. The generic characters are shown in Table \ref{generic character} and we have $\sum_{\gamma_z\in G^{+}_{1}}\chi^+(\gamma_z)=0$, whether $f=-N'$ or not, since $\left(\frac{2}{\cdot}\right)$ or $\left(\frac{-2}{\cdot}\right)$ is in the set of generic characters of discriminant $N'$.

\medskip
\paragraph*{\bf (3) $N\equiv 3(\text{mod }4)$.}
In this case, $-4N=N'l^2$ with $2\nmid N'$ and $2\mid l$. So the set of generic characters of $N'$ is $\left\{\left(\frac{\cdot}{p}\right)\Big|p\mid N'\right\}$. Recall that $G^{+}_{1}={(G^{+}_{1})}^{\prime}\bigsqcup{(G^{+}_{1})}^{\prime\prime}$ and there are two $1$-$1$ correspondences, between ${(G^{+}_{1})}^{\prime}$ and $\mathscr{Q}_{-4N}^0/\mathrm{SL}_2(\Z)$, and between ${(G^{+}_{1})}^{\prime\prime}$ and $\mathscr{Q}_{-N}^0/\mathrm{SL}_2(\Z)$.
Define
$$S^{\prime}_{f,+}=\left\{\left(
                     \begin{array}{cc}
                       a\sqrt{N} & -b/\sqrt{N} \\
                       c\sqrt{N} & -a\sqrt{N} \\
                     \end{array}
                   \right)
\in {(G^{+}_{1})}^{\prime}\Bigg|\left(\frac{b}{f}\right)=+1\right\}$$
and
$$S^{\prime\prime}_{f,+}=\left\{\left(
                     \begin{array}{cc}
                       a\sqrt{N} & -b/\sqrt{N} \\
                       c\sqrt{N} & -a\sqrt{N} \\
                     \end{array}
                   \right)
\in {(G^{+}_{1})}^{\prime\prime}\Bigg|\left(\frac{b/2}{f}\right)=+1\right\}.$$
Similarly, we can define $S^{\prime}_{f,-}$ and $S^{\prime\prime}_{f,-}$ by replacing the value $+1$ of Jacobi symbol by $-1$.
If $f\neq-N'$, the cardinality of $S^{\prime}_{f,+}$(resp. $S^{\prime\prime}_{f,+}$) is equal to that of $S^{\prime}_{f,-}$(resp. $S^{\prime\prime}_{f,-}$) and then $\sum_{\gamma_z\in G^{+}_{1}}\chi^+(\gamma_z)=0$.

If $f=-N'$, we have ${(G^{+}_{1})}^{\prime}=S^{\prime}_{f,+}$, ${(G^{+}_{1})}^{\prime\prime}=S^{\prime\prime}_{f,+}$,  and then
\begin{align*}
\sum_{\gamma_z\in {(G^{+}_{1})}^{\prime}}\chi^+(\gamma_z)&=\chi^+(W_N)h(-4N),\\
\sum_{\gamma_z\in {(G^{+}_{1})}^{\prime\prime}}\chi^+(\gamma_z)&=\chi^+(W_N)\chi(2)
\sum_{\gamma_z\in {(G^{+}_{1})}^{\prime\prime}}\chi\left(\frac{b}{2}\right)=\chi^+(W_N)\chi(2)h(-N).
\end{align*}
Therefore
$$\sum_{\gamma_z\in G^{+}_{1}}\chi^+(\gamma_z)=\sum_{\gamma_z\in {(G^{+}_{1})}^{\prime}}\chi^+(\gamma_z)+\sum_{\gamma_z\in {(G^{+}_{1})}^{\prime\prime}}\chi^+(\gamma_z)=\chi^+(W_N)\left(h(-4N)+\left(\frac{2}{f}\right)h(-N)\right)$$ as desired.
\medskip
\paragraph*{\bf (4) $N\equiv 0(\text{mod }4)$.}
For even $v_2(N)\geq2$, we have $-4N=N'l^2$ with $2\mid l$. Hence $\sum_{\gamma_z\in G_{1}^{+}}\chi^+(\gamma_z)=0$ if $f\neq-N'$ and $\chi^+(W_N)h(-4N)$ otherwise. If $v_2(N)\geq3$ is odd, then $-4N=N'l^2$ with $8\parallel N'$ and $2\mid l$, and then $\sum_{\gamma_z\in G_{1}^{+}}\chi^+(\gamma_z)=\chi^+(W_N)h(-4N)$ if $f=-N'$ and $\chi_2=\left(\frac{-2}{\cdot}\right)$, and $0$ otherwise.
\end{proof}

As an application, we relate $\text{dim }\mathcal{S}_k(\Gamma_0^+(N),\chi^+)$ with
$\text{dim }\mathcal{S}_k(\Gamma_0(N),\chi)$. Recall that
\[\nu_3(\chi)=\sum_{s\in A_3(N)}\chi(s),\ \nu_\infty(\chi)=\sum_{\substack{c|N \\\gcd\left(c,N/c\right)|N/f}}\varphi(\gcd(c,N/c)),\]
where $A_3(N)=\{s\in\mathbb{Z}/N\mathbb{Z}:s^2+s+1\equiv 0(\text{mod } N)\}$ (see Theorem $1$ of \cite{Cohen1976Dimension}).

\begin{Cor}\label{main2}
Let $N>3$ and $N\equiv 1(\text{mod }4)$. For quadratic character $\chi$ and integer $k\geq2$ with $\chi(-1)=(-1)^k$, $$\text{ dim }\mathcal{S}_k(\Gamma_0^+(N),\chi^+)
=\frac{1}{2}\text{ dim }\mathcal{S}_k(\Gamma_0(N),\chi).$$
\end{Cor}
\begin{proof}
By comparing dimension formulas \eqref{Gamma0N+equ} for $\Gamma^+_0(N)$ and that for $\Gamma_0(N)$ (see Theorem $1$ of \cite{Cohen1976Dimension}), it suffices to show
$$\nu_2^+(\chi^+)=\frac{\nu_2(\chi)}{2},\ \nu_3^+(\chi^+)=\frac{\nu_3(\chi)}{2},\ \nu_{\infty}^+(\chi^+)=\frac{\nu_{\infty}(\chi)}{2}.$$
The first term follows from Proposition \ref{helpproofcor1}. For the latter two terms, since the elliptic elements of $\Gamma_0(N)$ paired by $W_N$ have the same character value, then $\nu_3^+(\chi^+)=\frac{\nu_3(\chi)}{2}$.
The conclusion also holds for cusps and then $\nu_{\infty}^+(\chi^+)=\frac{\nu_{\infty}(\chi)}{2}$.
The statement then follows.
\end{proof}

\section{Dimension formulas of modular forms for $\Gamma_0^*(N)$ with character}\label{section 3}

In this section, we treat the group $\Gamma_0^*(N)$ and content ourselves with square-free $N$ from now on. We first prove the dimension formula for $\Gamma_0^*(N)$.

\begin{Thm}\label{main3}
Let $N>1$ be square-free, $k\in \mathbb{Z}$ with $\chi(-1)=(-1)^k$ and $\chi^*$ be a character of $\Gamma^*_0(N)$ extended from a trivial or quadratic Dirichlet character $\chi$ modulo $N$ for $\Gamma_0(N)$. We have
\begin{align}
  \text{ dim }&\mathcal{S}_k(\Gamma_0^*(N),\chi^*)-\text{ dim }\mathcal{M}_{2-k}(\Gamma_0^*(N),\overline{\chi^*})\label{Gamma0N*equ}\\
  &=\frac{k-1}{12}\nu_0^*(N)+\gamma_4(k)\nu_2^*(\chi^*)+
\gamma_3(k)\nu_3^*(\chi^*)+\gamma_8(k)\delta_8^*+\gamma_{12}(k)\delta^*_{12}-\frac{1}{2},\notag
\end{align}
where
\begin{align*}
\nu_0^*(N)&=2^{-\omega(N)}N\prod\limits_{p\mid N}\left(1+\frac{1}{p}\right),\\
\nu_2^*(\chi^*)&=\sum_{z\in E^*_{2}}\chi^*(\gamma_z),\\
\nu_3^*(\chi^*)&=\left\{
\begin{array}{lcl}
1 & & {\text{if } p\equiv 1\text{(mod }3) \text{ for 	all } p\mid N}\\
0 & & {\text{otherwise}}\\
\end{array} \right.,\\
\delta^*_8&=\left\{\begin{array}{lcl}
1 & & {\text{if } 2\mid N \text{ and } p\equiv 1 (\text{mod } 4)}\text{ for all }p\mid \frac{N}{2}\\
0 & & {\text{otherwise}}
\end{array} \right.,\\
\delta^*_{12}&=\left\{\begin{array}{lcl}
1 & & {\text{if } 3\mid N \text{ and } p\equiv 1 (\text{mod }3)}\text{ for all }p\mid \frac{N}{3}\\
0 & & {\text{otherwise}}
\end{array} \right.,
\end{align*}
$E^*_{e}$ is the set of $\Gamma^*_0(N)$-orbits of elliptic points of order $e$ and $\gamma_z$ is one generator of $\Gamma^*_0(N)_{\tau}/\{\pm I\}$ with $z=\pi(\tau)$.
\end{Thm}
%\subsection{Proof of Theorem \ref{main3}}
\begin{proof}
Substitute the signature of $\Gamma_0^*(N)$ (See Theorem 4 and Theorem 7 of \cite{Maclachlanternary1981}) to formula (3.3.6) in \cite{Shimura1974Onthetrace} and we get
\begin{align}
  \text{ dim }\mathcal{S}&_k(\Gamma_0^*(N),\chi^*)-\text{ dim }\mathcal{M}_{2-k}(\Gamma_0^*(N),\overline{\chi^*})\label{Shimuragamma0N*}\\
=&\frac{k-1}{12}\nu_0^*(N)-\frac{\nu_2^*(N)}{4}-\frac{\nu_3^*(N)}{3}-\frac{3\nu_4^*(N)}{8}-\frac{5\nu_6^*(N)}{12}-\frac{1}{2}+\sum_{z\in \mathfrak{R}^*_{N}}\mu_z',\notag
\end{align}
where $\mathfrak{R}^*_{N}$ consists of all elliptic points and cusps of $X(\Gamma_0^*(N))$, and $\nu^*_{e}(N)$ is the number of order-$e$ elliptic points with $e=2,3,4,6$.
Compare equations \eqref{Gamma0N*equ} and \eqref{Shimuragamma0N*} and then we need to prove
$$\sum\limits_{z\in E^*_{2}}\left(\mu_z'-\frac{1}{4}\right)=\gamma_4(k)\nu^*_2(\chi^*),\ \sum\limits_{z\in E^*_{3}}\left(\mu_z'-\frac{1}{3}\right)=\gamma_3(k)\nu^*_3(\chi^*),\ \mu'_{\infty}=0.$$
If there is order-$4$ or order-$6$ elliptic point, we have to prove additionally
$$\sum\limits_{z\in{E^*_{4}}}\left(\mu_z'-\frac{3}{8}\right)=\gamma_8(k)\delta^*_{8},\quad \sum\limits_{z\in{E^*_{6}}}\left(\mu_z'-\frac{5}{12}\right)=\gamma_{12}(k)\delta^*_{12}. $$
The proofs of these statements are analogous to those in the case of $\Gamma_0^+(N)$ which were given in the proof of Theorem \ref{main1} and we leave the details to the reader.
\end{proof}

For the rest of this section, let $N$ be square-free with $p\equiv1(\text{mod }4)$ for any prime factor $p\mid N$ and $e$ divides $N$ with $e>1$. We shall prove a $2$-power relation between the dimension of $\mathcal{S}_k(\Gamma_0^*(N),\chi^*)$ and that of $\mathcal{S}_k(\Gamma_0(N),\chi)$.

Set
$$\mathscr{Q}_{N,\Delta}^0=\left\{(a',b',c')\in \mathscr{Q}^0_{\Delta}|a'\equiv 0(\text{mod }N)\right\}$$
and recall that $\Gamma^e_0(N)$ is generated by $\Gamma_0(N)$ and $W_e$. Let $E^{e}_{2}$ be the set consisting of $\Gamma^e_0(N)$-orbits of elliptic points of order $2$ fixed by elements of $W_e\Gamma_0(N)$, whose cardinality is equal to $\prod_{p|N/e}\left(1+\left(\frac{-4e}{p}\right)\right)h(-4e)$ (see \cite{kluit}). Define $G^{e}$ the set of generators $\gamma_z$, one for each $z=\pi(\tau)\in E^{e}_{2}$, in $\Gamma^e_0(N)_\tau/\{\pm I\}$. Helling\cite{Helling1970} proved the 1-1 correspondence between $G^{e}$ and $\mathscr{Q}^0_{N,-4e}/\Gamma_0(N)$ given by
$\left(\begin{smallmatrix}
a\sqrt{e}&-b/\sqrt{e}\\cN/\sqrt{e}&-a\sqrt{e}\end{smallmatrix}\right)\mapsto (cN,-2ae,b).$
Note that $(cN,-2ae,b)$ is primitive, since $\gcd(b,ae)=1$ and $e\equiv1\text{(mod }4)$. Hence we can further require that $\gcd(b,-4e)=1$ for any element $\left(
   \begin{smallmatrix}
     a\sqrt{e} & -b/\sqrt{e} \\
     cN/\sqrt{e} & -a\sqrt{e} \\
   \end{smallmatrix}
 \right)\in G^{e}$, and then apply the properties of generic characters to determine the character sums on extra elliptic elements for $\Gamma^e_0(N)$.

Clearly $b'\text{ mod }2N$ is a $\Gamma_0(N)$-invariant of $(a',b',c')\in \mathscr{Q}^0_{N,-4e}$. Let $\varrho$ mod $2N$ satisfy $\Delta\equiv \varrho^2(\text{mod }4N)$ and set
$$\mathscr{Q}^0_{N,\Delta,\varrho}=\left\{(a',b',c')\in \mathscr{Q}^0_{\Delta}|a'\equiv 0(\text{mod }N),\ b'\equiv \varrho(\text{mod }2N)\right\},$$
$$G^{e}_{\varrho}=\left\{\left(
                     \begin{array}{cc}
                       a\sqrt{e} & -b/\sqrt{e} \\
                       cN/\sqrt{e} & -a\sqrt{e} \\
                     \end{array}
                   \right)
\in G^{e}\Bigg|-2ae\equiv\varrho(\text{mod }2N)\right\}.$$
\begin{Lem}\label{Gamma0Nesum}
Suppose $\chi$ is quadratic with conductor $N$, we have $\sum_{\gamma_z\in G^{e}}\chi^*(\gamma_z)=0$ for any $e\mid N$ with $e>1$.
\end{Lem}
\begin{proof}
Write $\gamma_z=W_e\gamma_z'\in G^{e}$ of the form $\left(\begin{smallmatrix}
                a\sqrt{e} & -b/\sqrt{e} \\
                cN/\sqrt{e} & -a\sqrt{e} \\
              \end{smallmatrix}
            \right)=\left(\begin{smallmatrix}
                a'\sqrt{e} & 1/\sqrt{e} \\
                c'N/\sqrt{e} & \sqrt{e} \\
              \end{smallmatrix}
            \right)\left(
              \begin{smallmatrix}
                a'' & b'' \\
                c''N & d'' \\
              \end{smallmatrix}
            \right)$ with $\gcd(b,-4e)=1$.
By the proposition of page $505$ in \cite{ZagierGross}, for fixed $\varrho$ mod $2N$ with $\varrho^2\equiv -4e(\text{mod }4N)$, there is an isomorphism $\mathscr{Q}^0_{N,-4e,\varrho}/\Gamma_0(N)\simeq\mathscr{Q}^0_{-4e}/\mathrm{SL}_2(\Z)$.
If we define
$$G^{e}_{\varrho,+}=\left\{\left(
                     \begin{array}{cc}
                       a\sqrt{e} & -b/\sqrt{e} \\
                       cN/\sqrt{e} & -a\sqrt{e} \\
                     \end{array}
                   \right)
\in G^{e}\Bigg|-2ae\equiv\varrho(\text{mod }2N),\ \left(\frac{b}{e}\right)=+1\right\},$$
$$G^{e}_{\varrho,-}=\left\{\left(
                     \begin{array}{cc}
                       a\sqrt{e} & -b/\sqrt{e} \\
                       cN/\sqrt{e} & -a\sqrt{e} \\
                     \end{array}
                   \right)
\in G^{e}\Bigg|-2ae\equiv\varrho(\text{mod }2N),\ \left(\frac{b}{e}\right)=-1\right\},$$
then $G^{e}_{\varrho}=G^{e}_{\varrho,+}\bigsqcup G^{e}_{\varrho,-}$, and the cardinality of $G^{e}_{\varrho,+}$ is equal to that of $G^{e}_{\varrho,-}$ by such isomorphism, Theorem 4.16 and Proposition 4.18 in \cite{binaryform}.
Hence
\begin{align*}
  \sum\limits_{\gamma_z\in G^{e}_{\varrho}}\chi^*(\gamma_z)=&\chi^*(W_e)\sum\limits_{\gamma_z\in G^{e}_{\varrho}}\chi(\gamma_z')
  =\chi^*(W_e)\sum\limits_{\gamma_z\in G^{e}_{\varrho}}\left(\frac{d''}{N}\right)\\
  =&\chi^*(W_e)\sum\limits_{\gamma_z\in G^{e}_{\varrho}}\prod_{p\mid e}\left(\frac{d''}{p}\right)\prod_{p\mid N/e}\left(\frac{d''}{p}\right)\\
  =&\chi^*(W_e)\sum\limits_{\gamma_z\in G^{e}_{\varrho}}\prod_{p\mid e}\left(\frac{a'b''e+d''}{p}\right)\prod_{p\mid N/e}\left(\frac{b''c'N/e+d''}{p}\right)\\
  =&\chi^*(W_e)\sum\limits_{\gamma_z\in G^{e}_{\varrho}}\left(\frac{b}{e}\right)\left(\frac{-a}{N/e}\right)\\
  =&\chi^*(W_e)\left(\frac{2e\varrho}{N/e}\right)\sum\limits_{\gamma_z\in G^{e}_{\varrho,+}\bigsqcup G^{e}_{\varrho,-}}\left(\frac{b}{e}\right)=0.
\end{align*}

Now we count the number of $\varrho$ mod $2N$ satisfying $\varrho^2\equiv -4e(\text{mod }4N)$ to compute $\sum\limits_{\gamma_z\in G^{e}}\chi^*(\gamma_z)$. 
Indeed, passing to prime factors $p$ of $4N$, we see that modulo $4$ there exists $1$ solution, and $\varrho^2\equiv -4e(\text{mod }p)$ has $1$ solution if $p\mid e$ and $\left(1+\left(\frac{-4e}{p}\right)\right)$ solutions if $p\mid \frac{N}{e}$. It follows that $$\sum\limits_{\gamma_z\in G^{e}}\chi^*(\gamma_z)=\prod\limits_{p\mid N/e}\left(1+\left(\frac{-4e}{p}\right)\right)\sum_{\gamma_z\in G^{e}_{\varrho}}\chi^*(\gamma_z)=0,$$
which finishes the proof.
\end{proof}
Put $E^*_{2}=E^{*}_{2,0}\bigsqcup E^{*}_{2,1}$, where $E^{*}_{2,0}$ and $E^{*}_{2,1}$ consist of $\Gamma^*_0(N)$-orbits of elliptic points of order $2$ fixed by elements of $\Gamma_0(N)$ and  $\Gamma^*_0(N)-\Gamma_0(N)$ respectively.
Let $G^{*}_{0}$ and $G^{*}_{1}$ be the set of generators $\gamma_z$ in $\Gamma_0^*(N)_\tau/\{\pm I\}$, one for each $z=\pi(\tau)\in E^{*}_{2,0}$ and $E^{*}_{2,1}$ respectively. Define  $G^*=G^{*}_{0}\bigsqcup G^{*}_{1}$.
\begin{Lem}\label{helpproofcoro2} If $\chi$ is quadratic with conductor $N$, then
$\sum_{\gamma_z\in G^{*}_{1}}\chi^*(\gamma_z)=0$.
\end{Lem}
\begin{proof}
For any $e\mid N$, we first prove that the set of order-$2$ elliptic points of $\Gamma_0^e(N)$ fixed by elements in $W_e\Gamma_0(N)\subset\Gamma_0^e(N)$ corresponds in $2^{\omega(N)-1}$-to-$1$ to the
the set of order-$2$ elliptic points of $\Gamma_0^*(N)$ fixed by elements in $W_e\Gamma_0(N)\subset\Gamma_0^*(N)$. Since the index $[\Gamma_0^*(N):\Gamma_0^e(N)]=2^{\omega(N)-1}$, this follows easily from the the Riemann-Hurwitz formula for the projection $X(\Gamma^e_0(N))\rightarrow X(\Gamma_0^*(N))$ and the following claim: Given any such point $z$ in the former and any $g\neq 1,e$ with $g\mid N$, the elliptic points $W_gz\not\sim z$ for $\Gamma_0^e(N)$. Suppose $W_gz\sim z$, so there exists $\gamma\in \Gamma_0^e(N)$ such that $\gamma W_gz=z$ with $\gamma W_g\in \Gamma_0^*(N)$ but $\gamma W_g\not\in\Gamma_0^e(N)$. Hence $z$ is an elliptic point for $\Gamma_0^*(N)$ of order at least that of the subgroup $\langle \gamma W_g, \Gamma_0^e(N)_z\rangle$ which contains at least $4$ elements, contradicting to the fact that $\Gamma_0^*(N)$ have only order-$2$ and order-$3$ elliptic points (see Theorem $4$ in \cite{Maclachlanternary1981}). 

Since for any elliptic element $\gamma\in W_e\Gamma_0(N)$ and any $g\mid N$,  $\chi^*(\gamma)=\chi^*(W_g\gamma W_g^{-1})$, applying the $2^{\omega(N)-1}$-to-$1$ correspondence above and Lemma \ref{Gamma0Nesum},
\[\sum\limits_{\gamma_z\in G^{*}_{1}}\chi^*(\gamma_z)=2^{1-\omega(N)}\sum\limits_{1<e\mid N}\sum\limits_{\gamma_z\in G^{e}}\chi^*(\gamma_z)=0.\]
This completes the proof.
\end{proof}
With the above conclusion, we can relate $\text{dim }\mathcal{S}_k(\Gamma_0^*(N),\chi^*)$ with
$\text{dim }\mathcal{S}_k(\Gamma_0(N),\chi)$.
\begin{Cor}\label{main4}
Let $N>1$ be square-free such that for all $p\mid N$, $p\equiv1(\text{mod }4)$. Suppose $\chi$ is quadratic with conductor $N$. For even integer $k\geq2$, $$\text{ dim }\mathcal{S}_k(\Gamma_0^*(N),\chi^*)
=\frac{1}{2^{\omega(N)}}\text{ dim }\mathcal{S}_k(\Gamma_0(N),\chi).$$
\end{Cor}
\begin{proof}
By comparing dimension formulas  \eqref{Gamma0N*equ} for $\Gamma^*_0(N)$ and that for $\Gamma_0(N)$ (see Theorem $1$ of \cite{Cohen1976Dimension}), we have to prove that
$$\nu_2^*(\chi^*)=\frac{\nu_2(\chi)}{2^{\omega(N)}},\ \nu_3^*(\chi^*)=\frac{\nu_3(\chi)}{2^{\omega(N)}},\ \nu_{\infty}^*(\chi^*)=\frac{\nu_{\infty}(\chi)}{2^{\omega(N)}}.$$
All the $\Gamma_0(N)$-inequivalent elliptic points $\tau$ of order $2$ are $\Gamma_0^*(N)$-equivalent via $W_g\tau\sim \tau$ for $g\mid N$ and the same is true for the order-$3$ elliptic points and for the cusps. So we have
\[\nu_{\infty}^*(\chi^*)=\frac{\nu_\infty(\chi)}{2^{\omega(N)}},\quad
\nu_3^*(\chi^*)=\frac{\nu_3(\chi)}{2^{\omega(N)}},\quad \sum\limits_{\gamma_z\in G_{0}^{*}}\chi^*(\gamma_z)=\frac{\nu_2(\chi)}{2^{\omega(N)}}.\]
Moreover,
$$\nu_2^*(\chi^*)=\sum\limits_{\gamma_z\in G^*}\chi^*(\gamma_z)=\sum\limits_{\gamma_z\in G_{0}^{*}}\chi^*(\gamma_z)+\sum\limits_{\gamma_z\in G_{1}^{*}}\chi^*(\gamma_z),$$
but by Lemma \ref{helpproofcoro2}, $\sum_{\gamma_z\in G^{*}_{1}}\chi^*(\gamma_z)=0$, so we are done with the proof.
\end{proof}

Note that for general square-free $N$, we can also compute $\nu^*_2(\chi^*)$ by similar computations as that of Proposition \ref{helpproofcor1}, but the resulting formula is complicated and we omit it.

We end this paper with the following example where $N=221$.
\begin{Exa}
Let us treat the case when $N=221=13\cdot17$, $k=6$,  $\chi=\left(\frac{\cdot}{221}\right)$, $\chi^+$ with $\chi^+(W_N)=-1$, and $\chi^*$ with $\chi^*(W_{13})=\chi^*(W_{17})=-1$.
Using MAGMA, we get that $$\nu_2(\chi)=\chi(21)+\chi(47)+\chi(174)+\chi(200)=-4, \ \text{dim }\mathcal{S}_k(\Gamma_0(N),\chi)=104.$$

The extra elliptic elements of $\Gamma_0^+(N)$ correspond to the quadratic forms in Table \ref{example221}, and we have $$\nu_2^+(\chi^+)=\frac{1}{2}\nu_2(\chi)+\chi^+(W_N)\sum_{\gamma_z\in R_{221,+}\bigsqcup R_{221,-}}\left(\frac{b}{221}\right)=\frac{1}{2}\nu_2(\chi)$$
with $\gamma_z=\left(\begin{smallmatrix}
 a\sqrt{N} & -b/\sqrt{N} \\c\sqrt{N} & -a\sqrt{N} \\\end{smallmatrix}\right)$ and $\gcd(b,-4N)=1$, which illustrates Proposition \ref{helpproofcor1} (see the meaning of $R_{f,+}$ and $R_{f,-}$ in the proof of Proposition \ref{helpproofcor1}). It follows that
$$\text{ dim }\mathcal{S}_k(\Gamma_0^+(N),\chi^+)=\frac{5}{12}\cdot\frac{252}{2}-
\frac{1}{4}(\chi(21)+\chi(47)+0)-1=52=\frac{1}{2}\text{ dim }\mathcal{S}_k(\Gamma_0(N),\chi),$$
which illustrates Corollary \ref{main2}.
\begin{table}
\centering
\caption{The quadratic forms associated to extra elliptic points of $\Gamma_0^+(N)$ }\label{example221}
%\vspace{0.3cm}
\setlength{\abovecaptionskip}{0pt}
\setlength{\belowcaptionskip}{0pt}
\begin{tabular}{|ccc|c|}
\hline
  $\left(\frac{-1}{\cdot}\right)$ & $\left(\frac{\cdot}{13}\right)$ & $\left(\frac{\cdot}{17}\right)$ &
  $(cN,-2aN,b)$ of $\Delta=-4N$ \\\hline
  $+$ & $+$ & $+$ &\begin{tabular}{c} $(N,0,1)$ $(30N,-46N,3897)$\\
  $(393N,-8N,9)$  $(9N,-8N,393)$
  \end{tabular} \\\hline
  $+$ & $-$ & $-$ &\begin{tabular}{c}   $(6N,-2N,37)$ $(177N,-4N,5)$\\$(6N,2N,37)$ $(5N,-4N,177)$ \end{tabular}\\\hline
  $-$ & $+$ & $-$ & \begin{tabular}{c}   $(74N,2N,3)$  $(10N,-6N,199)$\\$(74N,-2N,3)$ $(10N,6N,199)$ \end{tabular}\\\hline
  $-$ & $-$ & $+$ & \begin{tabular}{c}   $(2N,-2N,111)$ $(722N,-14N,15)$\\$(59N,-4N,15)$ $(15N,-4N,59)$ \end{tabular}\\
  \hline
\end{tabular}
\end{table}
\noindent

Denote $e_1=13$, $e_2=17$. The number of extra elliptic points fixed by $\Gamma^*_0(N)-\Gamma_0(N)$ for $\Gamma^*_0(N)$ (See Theorem 7 of \cite{Maclachlanternary1981}) is $$\frac{h(-4N)}{2}+h(-4e_1)+h(-4e_2)$$
and the corresponding classes of quadratic forms are listed in Table \ref{extrapointsforGamma_0^*(N)}, Table \ref{e_1} and Table \ref{e_2} respectively. We have
\begin{align*}
\nu_2^*(\chi^*)=&\frac{\nu_2(\chi)}{4}+2^{-1}\sum\limits_{1<e\mid N}\sum\limits_{\gamma_z\in G^{e}}\chi^*(\gamma_z)\\
=&\frac{\nu_2(\chi)}{4}+2^{-1}\left(2\sum_{\gamma_z\in G^{e_{1}}_{-16e_1}}\chi^*(\gamma_z)+2\sum_{\gamma_z\in G^{e_2}_{-8e_2}}\chi^*(\gamma_z)+\sum_{\gamma_z\in G^{N}_{0}}\chi^*(\gamma_z)\right)=\frac{\nu_2(\chi)}{4},
\end{align*}
and then
$$\text{ dim }\mathcal{S}_k(\Gamma_0^*(N),\chi^*)=\frac{5}{12}\cdot\frac{252}{4}-
\frac{1}{4}(\chi(21)+0)-1/2=26=\frac{1}{4}\text{ dim }\mathcal{S}_k(\Gamma_0(N),\chi),$$
which illustrates Corollary \ref{main4}.
\end{Exa}

\begin{table}
\centering
\caption{$h(-4N)/2$ quadratic forms associated to extra elliptic points for $\Gamma_0^*(N)$}\label{extrapointsforGamma_0^*(N)}
%\setlength{\abovecaptionskip}{2pt}
%\setlength{\belowcaptionskip}{2pt}
%\vspace{0.3cm}
\begin{tabular}{|ccc|c|}
\hline
  $\left(\frac{-1}{\cdot}\right)$ & $\left(\frac{\cdot}{e_1}\right)$ & $\left(\frac{\cdot}{e_2}\right)$ &
 $(cN,-2aN,b)$ of $\Delta=-4N$ \\\hline
  $+$ & $+$ & $+$ &\begin{tabular}{c} $(N,0,1)$
  $(393N,-8N,9)$
  \end{tabular} \\\hline
  $+$ & $-$ & $-$ &\begin{tabular}{c}   $(6N,-2N,37)$ $(6N,2N,37)$ \end{tabular}\\\hline
  $-$ & $+$ & $-$ & \begin{tabular}{c}   $(74N,2N,3)$  $(74N,-2N,3)$ \end{tabular}\\\hline
  $-$ & $-$ & $+$ & \begin{tabular}{c}   $(2N,-2N,111)$ $(59N,-4N,15)$ \end{tabular}\\
  \hline
\end{tabular}
\end{table}

%\vspace{0.3cm}

\begin{table}
\begin{minipage}[h]{0.49\linewidth}
    \centering 
 \caption{Quadratic forms of $-4e_1$ }\label{e_1}
\resizebox{7.1cm}{!}{
\begin{tabular}{|cc|c|}
 \hline
  $\left(\frac{-1}{\cdot}\right)$ & $\left(\frac{\cdot}{e_1}\right)$ &
 $(cN,-2ae_1,b)$, $\Delta=-52$ \\\hline
  $+$ & $+$ &$(49N,-16e_1,1)$\\\hline
  $-$ & $-$ & $(7N,-16e_1,7)$ \\
  \hline
  \end{tabular}}
 \end{minipage}
 \hfill
 \begin{minipage}{0.49\linewidth}
  \centering \caption{Quadratic forms of $-4e_2$ }\label{e_2}
    \resizebox{7.8cm}{!}{\begin{tabular}{|cc|c|}
   \hline
  $\left(\frac{-1}{\cdot}\right)$ & $\left(\frac{\cdot}{e_2}\right)$ &
  $(cN,-2ae_2,b)$, $\Delta=-68$ \\\hline
  $+$ & $+$ & \begin{tabular}{c}   $(21N,-8e_2,1)$ $(42N,34e_2,9)$ \end{tabular} \\\hline
  $-$ & $-$ &  \begin{tabular}{c}   $(3N,-8e_2,7)$ $(7N,-8e_2,3)$ \end{tabular} \\\hline
\end{tabular}}
\end{minipage}
\end{table}

%\subsection*{Acknowledgments}
%The first author was partially supported by a grant of National Natural Science Foundation of China (Grant No. 12271123). The authors are very grateful to the referee for pointing out to us the connection of the dimension formulas with Hecke operators and for many other valuable comments on the original version of the manuscript. 

\bibliographystyle{amsplain}
\bibliography{paper}

\end{document}